\documentclass[11pt]{amsart}
\usepackage{amsmath,amssymb,amsthm}
\usepackage[hidelinks]{hyperref}

\newtheorem{theorem}{Theorem}
\newtheorem{lemma}[theorem]{Lemma}
\newtheorem{remark}{Remark}
\newtheorem{proposition}[theorem]{Proposition}
\newtheorem{corollary}[theorem]{Corollary}

\newcommand{\HH}{\mathcal H_3}
\newcommand{\C}{\mathbb C}
\newcommand{\Q}{\mathbb Q}
\newcommand{\one}{\mathbf 1}
\newcommand{\id}{\mathrm{id}}
\DeclareMathOperator{\Hom}{Hom}

\title{A fusion category admitting no Hermitian structure}
\author{Hanshi Yang}
\address{Northwood High School}
\email{mainhanshiyangemail@gmail.com}
\date{\today}

\begin{document}

\begin{abstract}
A Galois conjugate of the Haagerup fusion category \(\mathcal H_3\) admits no Hermitian structure, even if the forms on its fusion spaces are allowed to be indefinite. The obstruction is visible in a single coefficient of the associator for three copies of \(\rho\), where a Hermitian structure would make this coefficient the complex conjugate of the corresponding coefficient of the inverse. This answers the first question posed by Etingof, Gelaki, Nikshych, and Ostrik in Remark~9.4.7 of \emph{Tensor Categories} in the negative. We also construct further counterexamples among near-group fusion categories that have nontrivial fusion multiplicities. 
\end{abstract}

\maketitle

\section{Introduction}
What prevents a complex fusion category from admitting a Hermitian structure? For a unitary structure sometimes one can rule out positivity by computing a categorical dimension, but this argument does not apply when indefinite forms are allowed. Etingof, Gelaki, Nikshych, and Ostrik raised the question of whether a Hermitian structure can fail to exist at all \cite[Remark~9.4.7]{EGNO}.

Consider the six-object Haagerup category \(\HH\) of Grossman and Snyder \cite{GS}. We use the algebraic presentation of its associator given by Huang and Lin \cite[Section~5.2.2]{HL}. Let \(\overline{\Q}\subset\C\); let \(K\) be the number field generated by \(\sqrt{13}\) and the structure constants of this presentation, notably including a choice of algebraic duality maps. 
For an embedding \(\sigma:K\hookrightarrow\C\), let \(\HH^\sigma\) denote the category obtained by applying \(\sigma\) to these constants and extending the scalars to \(\C\). The same fusion rules will define a complex fusion category with the conjugated associator since under \(\sigma\), the pentagon and the duality identities are still preserved.

\begin{theorem}\label{thm:main}
For every embedding \(\sigma:K\hookrightarrow\C\) with
\(\sigma(\sqrt{13})=-\sqrt{13}\), the category \(\HH^\sigma\) admits
no Hermitian structure.
\end{theorem}

Such embeddings exist since the nontrivial embedding of
\(\Q(\sqrt{13})\) extends to \(K\). The theorem
addresses the first question in \cite[Remark~9.4.7]{EGNO}.

We use the dagger formulation of a Hermitian structure in which the associator and unit constraints are unitary.

For a simple object \(D\) and maps \(f,g:D\to W\), let
\[
    f^\dagger g=\langle f,g\rangle\,\id_D .
\]
The resulting Hermitian form is linear in the second argument and need not be positive definite.

Let \(X,D\) be simple objects in a multiplicity-free fusion category and choose a nonzero inclusion for each simple summand in each product of two simples. When one composes these inclusions, one gets tree bases \(L_e\) and \(R_f\) for maps from \(D\) to \((X\otimes X)\otimes X\) and to \(X\otimes(X\otimes X)\) respectively. Here \(e\) and \(f\) denote the intermediate simple objects. With \(\Phi\) the associator, the matrix
convention is
\[
    \Phi L_e=\sum_f R_f F_{f,e}.
\]

\begin{lemma}\label{lem:entry}
Suppose that \(X\) and \(D\) both occur in \(X\otimes X\). If the category admits a Hermitian structure, then
\[
     F_{X,X}=\overline{(F^{-1})_{X,X}} .
\]
\end{lemma}

\begin{proof}
Choose inclusions denoted \(u:X\to X\otimes X\) and \(v:D\to X\otimes X\). 

Let \(u^\dagger u=s\,\id_X\) and \(v^\dagger v=t\,\id_D\). Both scalars are real. To see that \(s\ne0\), choose a retraction \(p\) of \(u\). Because the space \(\Hom(X\otimes X,X)\) is one-dimensional, the nonzero map \(u^\dagger\) is a nonzero multiple of \(p\). The same argument gives \(t\ne0\). It follows that
\[
     \langle L_X,L_X\rangle=st=\langle R_X,R_X\rangle .
\]
Both tree bases are orthogonal since the intermediate simple objects are pairwise nonisomorphic.

Expanding \(\Phi L_X\) in the right tree basis gives \(st\,F_{X,X}\). Moving
\(\Phi\) to the first slot gives
\(\langle\Phi^{-1}R_X,L_X\rangle\). Expanding in the left tree basis
gives \(st\,\overline{(F^{-1})_{X,X}}\). Cancelling \(st\) proves the
identity.

Rescaling \(u\) and \(v\) by \(\lambda\) and \(\mu\), respectively, will rescale both \(L_X\) and \(R_X\) by \(\lambda\mu\). Hence both diagonal entries are unchanged.
\end{proof}

\section{Three copies of \texorpdfstring{\(\rho\)}{rho}}

The invertible simples of \(\HH\) are \(\one,\alpha,\alpha^2\), and the remaining simples are \(\rho_j=\alpha^j\rho\), with \(j\in\mathbb Z/3\mathbb Z\). The fusion rules are 
\[
    \alpha^3=\one,\qquad \rho\alpha=\alpha^{-1}\rho,\qquad
    \rho^2=\one\oplus\rho_0\oplus\rho_1\oplus\rho_2.
\]
Take  \(X=\rho\) and \(D=\rho_1\). Since the unit channel contributes only to \(\rho\), the intermediate simples in both tree bases are \(\rho_0,\rho_1,\rho_2\), which we order accordingly.

For \(D=\rho_1\), the entries in Huang and Lin's table with the first channel \(\rho_0\) are \(y_1,y_2,z\) \cite[Theorem~5.2(b), Section~5.2.2]{HL}. Let \(r=\sqrt{13}>0\) and \(q=\sqrt{6(1+r)}>0\). Then
\[
    y_1=\frac{5-r-q}{12},\qquad
    y_2=\frac{5-r+q}{12},\qquad z=\frac{1+r}{6}.
\]
By the channel-shift symmetry of \cite[Section~5.2, footnote~6]{HL}, the entry will be unchanged under
\[
    (\rho_i,\rho_j)\mapsto(\alpha\rho_i,\rho_j\alpha)
    =(\rho_{i+1},\rho_{j-1}).
\]
Hence the entries depend only on \(i+j\) modulo three and
\[
    F_0=
    \begin{pmatrix}
        y_1&y_2&z\\ y_2&z&y_1\\ z&y_1&y_2
    \end{pmatrix}.
\]

\begin{proof}[Proof of Theorem~\ref{thm:main}]
Let \(a=\sigma(y_1)\), \(b=\sigma(y_2)\), and \(c=\sigma(z)\).
Since \(\sigma(q)^2=6(1-r)<0\),
\[
  \{a,b\}=
  \left\{\frac{5+r+i\sqrt{6(r-1)}}{12},
         \frac{5+r-i\sqrt{6(r-1)}}{12}\right\},
  \qquad c=\frac{1-r}{6}.
\]
In particular \(a\) is nonreal for either extension on \(q\). The conjugated matrix \(F=\sigma(F_0)\) satisfies \(F^2=I\) and its row sums are \(a+b+c=1\), and
\[
    ab=\frac{(5+r)^2+6(r-1)}{144}
       =\frac{2+r}{9}=-c(a+b).
\]
Each off-diagonal entry of \(F^2\) is therefore \(ab+bc+ca=0\) and each diagonal entry is \(a^2+b^2+c^2=(a+b+c)^2=1\).

It follows that the coefficient from \(\rho\) to \(\rho\) is \(a\) in both the associator and its inverse. Lemma~\ref{lem:entry} would then give \(a=\overline a\), a contradiction.
\end{proof}

\begin{remark}
The conclusion of Theorem~\ref{thm:main} is already implicit in the work of Evans and Gannon \cite{EG}, although it appears to contradict the assertion of Hermitianity in \cite[Theorem~3(ii)]{EG}. In fact, the matrix \(A\) corresponding to the last row of their table has nonreal diagonal entries which, by their criterion in Theorem~2(c), obstruct Hermitianity.
\end{remark}

\section{Near-group examples}

The obstruction is not restricted to multiplicity-free fusion categories as we will now show. Consider the near-group categories of type \(G+n\) of Evans and Gannon \cite{EG2}, where \(G=\mathbb Z_n\). There is one noninvertible simple object \(\rho\), with fusion rules
\[
    g\rho=\rho g=\rho,\qquad
    \rho^2=\bigoplus_{g\in G}g\oplus n\rho.
\]
Use their endomorphism presentation with parameters \(a(g),b(g),c\). Let
\[
    d_\pm=\frac{n\pm\sqrt{n^2+4n}}{2}.
\]

Let \(\sigma\) be an embedding of the field of definition sending \(d_+\) to \(d_-\).

\begin{proposition}\label{prop:near-group}

If the conjugated category admits a Hermitian structure, then
\[
    |1-\sigma(\sqrt n\,c)|^2
    \leq\frac{(n-1)^2d_-^2}{n}.
\]
\end{proposition}

\begin{proof}
Suppose the conjugated category admits a Hermitian structure. Let \(T_g\) and \(T_g^\vee\) be the transported bases of \(\Hom(\rho,\rho^2)\) and \(\Hom(\rho^2,\rho)\), respectively, with \(T_g^\vee T_h=\delta_{g,h}\id_\rho\). Under the action of \(G\) on the \(T_g\) in \cite[Theorem~1]{EG2}, the \(T_g\) belong to distinct character spaces and must therefore be orthogonal. It is clear that the action of \(G\) preserves the Hermitian form since tensoring with the identity of an invertible object will commute with the dagger. Let
\[
    T_g^\dagger=t_g T_g^\vee,\qquad t_g\in\mathbb R^\times.
\]

Let \(U_g:\rho\to\rho\alpha_g\) be the morphism in \cite{EG2}. Tensoring \(U_g\) with \(\id_\rho\) gives an endomorphism \(M_g\) of \(\rho^2\) with \(M_gT_h=T_{h-g}\). Since \(\rho\) and \(\rho\alpha_g\) are simple, \(M_g^\dagger M_g\) will be a nonzero real scalar multiple of the identity. But \(M_g^n=\id_{\rho^2}\), so this scalar is either \(1\) or \(-1\). Hence \(t_g/t_0\in\{1,-1\}\).

Consider the two trees
\[
    L=(T_0\otimes\id_\rho)T_0,\qquad
    R_g=(\id_\rho\otimes T_g)T_0.
\]
By \cite[(2.6), Corollary~5]{EG2},
\[
    L^\vee R_g=\sigma(b(g))\id_\rho,\qquad
    R_g^\vee L=\sigma(\overline{b(g)})\id_\rho.
\]
Since \(L^\dagger=t_0^2L^\vee\) and \(R_g^\dagger=t_0t_gR_g^\vee\), 
\[
    t_0^2\overline{\sigma(b(g))}
    =t_0t_g\sigma(\overline{b(g)}).
\]
For \(g\ne0\) we have \(b(g)\overline{b(g)}=1/n\). Hence,
\[
    |\sigma(b(g))|^2=\frac{t_g}{nt_0}.
\]
The left-hand side is positive, so \(t_g=t_0\) and \(|\sigma(b(g))|=1/\sqrt n\).

Evaluating the corresponding identity in \cite{EG2} at the identity element gives
\[
    \sum_{g\ne0}b(g)=\frac{1-\sqrt n\,c}{d_+}.
\]
Applying \(\sigma\) and using the triangle inequality gives the necessary condition 
\[
    |1-\sigma(\sqrt n\,c)|^2
    \leq\frac{(n-1)^2d_-^2}{n}.
\]
\end{proof}

\begin{corollary}
For suitable embeddings sending \(d_+\) to \(d_-\), the same obstruction gives some more counterexamples from \cite[Table~2]{EG2}. Some examples are:
\begin{itemize}
    \item \(\mathcal C(\mathbb Z_3,3)\), first and second
    entries, with \(\sigma(\sqrt{21})=-\sqrt{21}\).

    \item \(\mathcal C(\mathbb Z_4,4)\), first entry,
    with \(\sigma(\sqrt2)=-\sqrt2\).

    \item \(\mathcal C(\mathbb Z_5,5)\), second and third
    entries, with \(\sigma(\sqrt5)=-\sqrt5\).

    \item \(\mathcal C(\mathbb Z_7,7)\), first entry,
    with \(\sigma(\sqrt{77})=-\sqrt{77}\).
\end{itemize}
\end{corollary}

\begin{proof}
Take the second \(\mathbb Z_3\) entry in \cite[Table~2]{EG2}, where \(c=e^{2\pi i/12}\), and choose an embedding \(\sigma\) fixing \(\sqrt3\) and sending \(\sqrt{21}\) to \(-\sqrt{21}\) and \(c\) to \(c^{-1}\). The difference between the two sides of the inequality is
\[
    1-\frac{4}{3}\left(\frac{3-\sqrt{21}}{2}\right)^2
    =2\sqrt{21}-9>0.
\]
Thus the conjugated category, which has four simple objects, admits no Hermitian structure.

The remaining cases follow by the same calculation.
\end{proof}

\begin{remark}
The examples above suggest that the obstruction to a Hermitian structure can be studied further. It would be interesting to determine whether these obstructions, or a suitable generalization of them, give a necessary and sufficient condition for the existence of a Hermitian structure. We leave this as an open remark for now. 
\end{remark}
 
\section*{Acknowledgments}
This work grew out of an earlier project, and I thank Du Pei for his mentorship. I also thank Pavel Etingof and David Reutter for correspondence concerning the obstruction and its relation to the work of Evans and Gannon. The presentation of the two-trees argument through \(\rho\) draws in part on an explanation that Etingof prepared with the assistance of ChatGPT.

\end{document}